\documentclass[12pt]{article}

\def\rcs $#1: #2 ${\expandafter\def\csname rcs#1\endcsname {#2}}
\rcs $Date: 2008/10/24 02:14:39 $

\usepackage{amsmath,amssymb,amsfonts,amsthm,url,xspace,graphicx}
\usepackage[pdfauthor={Omer Angel, Abraham D. Flaxman, David B. Wilson},
            pdftitle={A sharp threshold for minimum bounded-depth/diameter spanning and Steiner trees},
            bookmarks=true,bookmarksopen=true,bookmarksopenlevel=2]{hyperref}

\newcommand{\sref}[1]{\hyperref[#1]{\S~\ref*{#1}}}
\newcommand{\fref}[1]{\hyperref[#1]{Figure~\ref*{#1}}}
\newcommand{\tref}[1]{\hyperref[#1]{Theorem~\ref*{#1}}}
\newcommand{\lref}[1]{\hyperref[#1]{Lemma~\ref*{#1}}}
\newcommand{\cref}[1]{\hyperref[#1]{Corollary~\ref*{#1}}}
\newcommand{\pref}[1]{\hyperref[#1]{Proposition~\ref*{#1}}}

\newcommand{\lemref}[1]{\lref{L:#1}}

\newtheorem{theorem}{Theorem}[section]
  \newtheorem{lemma}[theorem]{Lemma}
  
  \newtheorem{corollary}[theorem]{Corollary}

\theoremstyle{definition}
  \newtheorem{remark}[theorem]{Remark}
  \newtheorem{definition}[theorem]{Definition}

\renewcommand{\th}{\ensuremath{^{\rm th}}}

\newcommand{\E}{\mathbb{E}}

\newcommand{\R}{\mathbb{R}}
\newcommand{\N}{\mathbb{N}}

\newcommand{\T}{\mathcal{T}}
\newcommand{\eps}{\varepsilon}
\newcommand{\Exp}{\operatorname{Exp}}

\newcommand{\old}[1]{}
\renewcommand{\th}{\ensuremath{^{\text{th}}}\xspace}

\newcommand{\1}{\mathbf{1}}

\newcommand{\Rdn}{R_{\delta,n}}
\newcommand{\fR}{f^{(\Rdn)}}
\newcommand{\eqd}{\stackrel{\mathcal D}{=}}

\newcommand{\wt}{\operatorname{wt}}
\newcommand{\MST}{\operatornamewithlimits{MST}}
\newcommand{\MSTrr}[2]{\underset{\text{\rm depth$(#1)\leq #2$}}{\MST}}
\newcommand{\MSTr}[1]{\underset{\text{\rm depth$\leq #1$}}{\MST}}
\newcommand{\MSTd}[1]{\underset{\text{\rm diam$\leq #1$}}{\MST}}
\begin{document}

\title{\vspace*{-80pt}
              A sharp threshold for \\
     minimum bounded-depth and bounded-diameter \\
 spanning trees and Steiner trees in random networks
}

\author{Omer Angel\thanks{\noindent University of Toronto.}
                  \thanks{University of British Columbia.}
 \and Abraham D. Flaxman\thanks{\noindent Microsoft Research.}
                        \thanks{University of Washington.}
 \and David B. Wilson\footnotemark[3]
}

\date{}

\maketitle

\vspace*{-18pt}
\begin{abstract}
  In the complete graph on $n$ vertices, when each edge has a weight
  which is an exponential random variable, Frieze proved that the
  minimum spanning tree has weight tending to
  $\zeta(3)=1/1^3+1/2^3+1/3^3+\cdots$ as $n\to\infty$.  We consider
  spanning trees constrained to have depth bounded by $k$ from a
  specified root.  We prove that if $k\geq\log_2\log n+\omega(1)$,
  where $\omega(1)$ is any function going to $\infty$ with $n$,
  then the minimum bounded-depth spanning tree still has weight
  tending to $\zeta(3)$ as $n\to\infty$, and that if $k<\log_2\log n$,
  then the weight is doubly-exponentially large in $\log_2\log n \,-\,k$.
  It is NP-hard to find
  the minimum bounded-depth spanning tree, but when
  $k\leq\log_2\log n -\omega(1)$, a simple greedy algorithm is
  asymptotically optimal, and when $k\geq\log_2\log n+\omega(1)$, an
  algorithm which makes small changes to the minimum (unbounded depth)
  spanning tree is asymptotically optimal.  We prove similar results
  for minimum bounded-depth Steiner trees, where the tree must connect
  a specified set of $m$ vertices, and may or may not include other
  vertices.  In particular, when $m={\rm const}\times n$, if
  $k\geq\log_2\log n+\omega(1)$, the minimum bounded-depth Steiner
  tree on the complete graph has asymptotically the same weight as the
  minimum Steiner tree, and if $1\leq k\leq\log_2\log n-\omega(1)$,
  the weight tends to $(1-2^{-k})\sqrt{8m/n}
  \big[\sqrt{2mn}/2^k\big]^{1/(2^k-1)}$ in both expectation and
  probability.  The same results hold for minimum bounded-diameter
  Steiner trees when the diameter bound is $2k$; when the diameter
  bound is increased from $2k$ to $2k+1$, the minimum Steiner tree
  weight is reduced by a factor of $2^{1/(2^k-1)}$.
\end{abstract}

\section{Introduction}

For graphs with random edge weights, we study minimum spanning trees
and Steiner trees in which there is a bound on the diameter or else a
bound on the depth from a specified root vertex.  We obtain precise
estimates of the weight of the optimal tree, and some of the results
are surprising.  There is a sharp cutoff in the depth/diameter constraint,
above which the constrained minimum spanning tree has almost the same
weight as the unconstrained minimum spanning tree, and below which the
weight blows up.

\subsection{Definitions}

The minimum spanning tree (MST) of an edge-weighted undirected simple
graph~$G$ on $n$ vertices is the spanning tree which minimizes the sum
of the edge weights.  The Steiner tree problem also specifies a
set~$T$ of $m=|T|\leq n$ terminal vertices that are to be connected by
the tree; the tree may or may not contain the other vertices in the
graph~$G$.  We denote the \underline{m}inimum \underline{s}panning
\underline{t}ree of~$G$ by $\MST(G)$, and the \underline{m}inimum
\underline{S}teiner \underline{t}ree by $\MST(G,T)$.  Prim's algorithm
and Kruskal's algorithm are two classic efficient algorithms for
finding the $\MST(G)$, but computing the minimum Steiner tree
$\MST(G,T)$ is well-known to be NP-hard \cite{Garey-Johnson}.

The bounded-depth Steiner tree problem, also known as the Steiner tree
problem with ``hop constraints,'' is an abstraction of several
important real-world combinatorial optimization problems, including
designing telecommunications networks with a maximum transmission
delay bound \cite{GouveiaM2003} and solving lot sizing problems with a
limits on the number of periods goods can in stock \cite{Voss-1999}.
We denote by $\MSTrr{r}{k}(G,T)$ the minimum weight Steiner tree of
the graph~$G$ connecting the terminal vertices~$T$ such that each
vertex is within distance~$k$ from the root vertex~$r$.  Similarly,
$\MSTd{k}(G,T)$ denotes the minimum weight Steiner tree of~$G$
connecting the vertices in~$T$ with diameter bounded by~$k$.  The
minimum bounded-depth and bounded-diameter spanning trees are of
course the special case where $T$ is the set of all vertices.

For a tree $\cal T$, we let $\wt(\cal T)$ denote its weight.

There has been extensive research (which we describe below) in
computer science, mathematics, operations research, and physics on the
bounded-diameter and bounded-depth versions of these problems.  In the
bounded-diameter version, the minimization is only over trees which
satisfy a bound on their diameter (maximum number of edges within the
tree connecting a pair of vertices), and in the bounded-depth
version, the minimization is over trees satisfying a bound on the
maximum distance from a pre-specified root vertex.  (The
bounded-diameter and bounded-depth versions are closely related.)

\subsection{The MST on random graphs}

There has been a lot of research on the properties of MST's on random
graphs.  Two of the most well-studied ensembles of random graphs are
the following:
\begin{enumerate}
\item The vertices of $G$ are the points of a Poisson point process in
  Euclidean space, with the edge weights being the Euclidean distance
  between the points.  The MST for these geometric graphs~$G$ has been
  studied in \cite{Avram-Bertsimas} \cite[Chapter~13]{Penrose}
  \cite{ClementiILMRS2007} and other articles.
\item The graph~$G$ is the complete graph~$K_n$, with the edge weights
  being i.i.d.\ copies of a random variable, such as an exponential with
  mean~$1$, or a uniform number between $0$ and~$1$.  (It turns out to
  matter very little which random variate occurs on
 the edges.)
\end{enumerate}

In 1985, Frieze showed that the expected cost of the minimum spanning
tree on the complete graph with edge weights distributed independently
and uniformly between 0 and~1 tends to a constant as $n$ tends to~$\infty$,
and the constant is
$\zeta(3)=1/1^3+1/2^3+1/3^3+\cdots=1.202\ldots$ \cite{Frieze1985}.  In our
notation, this says $\E[\wt(\MST(K_n))] \rightarrow \zeta(3)$ as $n
\rightarrow \infty$ (we abuse notation by making the edge weight
distribution implicit in the graph $K_n$).  A concentration result was
also proven, so the actual weight is with high probability close to
$\zeta(3)$ \cite{Frieze1985}.  Even more precise results are known:
the distribution of $\wt(\MST(K_n))$ converges to a Gaussian with mean
$\zeta(3)$ and variance $(6\zeta(4)-4\zeta(3))/n$
\cite{MR1369071,MR2225705}.  Since most edges in the optimal tree have
weight close to~$0$, so long as the weight random variables have a
density function that is~$1$ at weight~$0$, $\wt(\MST(K_n)) \to
\zeta(3)$ with high probability \cite{MR905183,MR1054012}.

\nocite{MR2007286}

Regarding the structure of $\MST(K_n)$, it is known that, with high
probability, the diameter of $\MST(K_n)$ is $\Theta(n^{1/3})$, and the
expected diameter is also $\Theta(n^{1/3})$ \cite{ABBR}.  (This is in
contrast to the uniformly random spanning tree on the complete graph,
which has diameter $\Theta(n^{1/2})$ in probability and in expectation
\cite{renyi-szekeres,szekeres}.)  From this diameter bound, it follows
that $\wt\!\Big(\MSTrr{r}{\omega(n^{1/3})}(K_n)\Big) \rightarrow \zeta(3)$
and
$\wt\!\Big(\MSTd{\omega(n^{1/3})}(K_n)\Big) \rightarrow \zeta(3)$
in probability.  (When $G$ is the complete graph $K_n$, we can
drop the root vertex $r$ from the notation.)  We prove that
this convergence still holds for a much more restrictive diameter
bound or depth bound:
\begin{theorem} \label{MST-cutoff}
For the complete graph~$K_n$ with $\Exp(1)$ edge weights,
if $k=\log_2\log n + \omega(1)$, where $\omega(1)$ is any quantity that
tends to~$\infty$, however slowly, then
\begin{equation}
\wt\!\Big(\MSTr{k}(K_n)\Big) \rightarrow
\zeta(3)\ \ \text{and}\ \ \wt\!\Big(\MSTd{2k}(K_n)\Big) \rightarrow
\zeta(3)
\end{equation}
in both probability and expectation.
This is tight in the sense that when $k=\log_2\log n - \Delta$,
\begin{equation}
\wt\!\Big(\MSTr{k}(K_n)\Big) = \exp(2^{\Delta+\Theta(1)})\ \ \text{and}\ \ 
\wt\!\Big(\MSTd{2k}(K_n)\Big) = \exp(2^{\Delta+\Theta(1)}).
\end{equation}
\end{theorem}
Thus there is a sharp cutoff at depth $\log_2\log n \pm \Theta(1)$
(diameter $2\log_2\log n \pm \Theta(1)$).

\begin{remark}
Any tree of depth~$k$ has diameter $\leq 2k$.  Any tree of diameter
$2k$ is also a tree of depth~$k$, rooted at a uniquely defined central
vertex.  Similarly, any tree of diameter $2k+1$ is a tree of depth~$k$
rooted at a central edge.  Thus, the principal difference between
$\MSTr{k}$ and $\MSTd{2k}$ is that in the first case, the root vertex
is pre-specified, and in the second case, any vertex may serve as the
root.  Thus, for a given set $T$ of terminal vertices,
\begin{equation}\label{diam<depth}
 \wt\!\Big(\MSTd{2k}(G,T)\Big)\leq\wt\!\Big(\MSTrr{r}{k}(G,T)\Big).
\end{equation}
There are some parameter values for which the bounded-diameter
tree is quite a bit lighter than the bounded-depth tree.  For example,
if $k=1$ and $T$ consists of $2$ vertices $u$ and~$v$, we have
$\wt\!\Big(\MSTd{2}(K_n,\{u,v\})\Big)=\Theta(1/\sqrt{n})$ while
$\wt\!\Big(\MSTr{1}(K_n,\{u,v\})\Big)=\Theta(1)$.  But for the parameter
values covered by our theorems, it turns out that
$\wt\!\Big(\MSTd{2k}(K_n,T)\Big)=(1+o(1))\wt\!\Big(\MSTr{k}(K_n,T)\Big)$.
\end{remark}

\begin{remark}
  This theorem also holds when the edge weights come from other
  distributions, as discussed below.
\end{remark}

\begin{remark}
  Since our focus is on the complete graph $K_n$, we can simplify
  notation for the minimum Steiner tree by writing $\MST(K_n,m)$
  instead of $\MST(K_n,T)$ for a set~$T$ of $m$ terminal vertices.
\end{remark}

Steiner trees in networks with uniformly random edges weights on the
complete graph~$K_n$ were
investigated in \cite{BGRS}.  Since $K_n$ is symmetric, it is only necessary
to specify the number~$m$ of terminals rather than the precise set.
There it was shown that when $2\leq m\leq o(n)$,
$$\E[\wt(\MST(K_n,m))] = (1+o(1))\frac{m-1}{n} \log \frac{n}{m}.$$
When $m$ is of the same order as $n$, say
$m=\alpha n$, the value of $\E[\wt(\MST(K_n,\alpha n))]$ is
$\Theta(1)$, but its precise limiting value is not known except when
$\alpha=1$ (where it is $\zeta(3)$).  We prove that this same weight
can be achieved when suitably restricting the diameter:
\begin{theorem} \label{Steiner-deep}
  For any number $m$ of terminal vertices ($2\leq m\leq n$),
  if $k\geq \log_2\log n + \omega(n/(m \log (e n/m)))$,
  and the edge weight probability distribution has density~$1$ at~$0$, then
\begin{equation}
 \frac{\wt\!\Big(\MSTr{k}(K_n,m)\Big)}{\wt(\MST(K_n,m))} \to 1
 \ \ \text{and} \ \ 
 \frac{\wt\!\Big(\MSTd{2k}(K_n,m)\Big)}{\wt(\MST(K_n,m))} \to 1
\end{equation}
in probability, and if the expected edge weight is finite, convergence
holds in expectation too.
\end{theorem}
In the case $m=n$, \tref{Steiner-deep} yields the first part of
\tref{MST-cutoff}, but for Steiner trees with general $m$ we do
not know if there is a sharp cutoff in the same sense that there is
for the minimum spanning tree, though whenever $m=\Theta(n)$ there is
still a sharp cutoff at $\log_2\log n \pm\Theta(1)$.

The next theorem gives the weight when the depth is smaller than
$\log_2\log n$.

\begin{theorem}\label{th:k-small}
  If there are $m=n^{1-o(1)}$ terminal vertices,
  and $2\leq k<\log_2\log n-\log_2\log(e n/m)-\omega(1)$,
  and the edge weight probability distribution has density~$1$ at~$0$, then
\begin{align}
\left.
\begin{array}{r}
\wt\!\Big(\MSTr{k}(K_n,m)\Big) \\
\wt\!\Big(\MSTd{2k}(K_n,m)\Big)
\end{array}
\right\}
& = (1-2^{-k} \pm o(1)) \sqrt{\frac{8m}{n}}
\left( \frac{\sqrt{2mn}}{2^{k}} \right)
    ^{\frac{1}{2^k-1}} \\
\wt\!\Big(\MSTd{2k+1}(K_n,m)\Big)
&= (1-2^{-k} \pm o(1)) \sqrt{\frac{8m}{n}}
\left( \frac{\sqrt{mn/2}}{2^{k}} \right)
    ^{\frac{1}{2^k-1}}
\end{align}
  in probability, and if the expected edge weight is finite, convergence
  holds in expectation too.  (These formulas are valid for $k=1$ too,
  when the weights are $\Exp(1)$ random variables.)
\end{theorem}

For example, when $k=2$ and $m=\alpha n$, \tref{th:k-small}
implies that with high probability
$$\wt\!\Big(\MSTr{2}(K_n,\alpha n)\Big) = (1+o(1)) \frac32 \alpha^{2/3} n^{1/3}.$$

The second part of \tref{MST-cutoff} follows from \tref{th:k-small}
upon specializing to the case $m=n$ and $k=\log_2\log n - \Delta$.

\subsection{Computational intractability and approximation algorithms}

Obtaining optimal trees is computationally intractable in general.
The minimum bounded-diameter spanning tree problem is NP-hard for
any diameter between $4$ and $n-2$ \cite[pg.~206]{Garey-Johnson}, and
the minimum bounded-depth spanning tree problem is NP-hard even
for depth~$2$ (this can be shown by a reduction from the facility
location problem, see \cite{DGR}).  Of course the bounded-diameter
and bounded-depth Steiner tree problems are only harder, so they too
are NP-hard.  In fact, for any fixed diameter $\geq 4$, it is
NP-hard to even approximate the minimum bounded-diameter spanning
tree to within an approximation ratio of better than $O(\log n)$
\cite{MR1795243}.

Because of this complexity, numerous algorithms have been
investigated, including exact (but time-consuming) integer programming
formulations \cite{MR1165807,MR1240161,GR05a}, fast rigorous
approximation algorithms \cite{MR1795243,MR1700182,MR2101043}, and
heuristic approximation algorithms
\cite{DGR,Voss-1999,Gouveia95,Monnot,Gouveia,CCL,CCLb,RJ03a,RJ03b,
GR05b,GHR06,kopinitsch,putz,zaubzer,BayatiBBCRZ2008}.

These intractability and inapproximability results are of course for
worst-case graphs, and for random graphs one can do better.  One of
the heuristic algorithms, based on ``survey propagation'' and ``the
cavity method,'' was recently tested on random graphs
\cite{BayatiBBCRZ2008}, which led us to investigate the weight of the
true optimal minimum bounded-depth and bounded-diameter spanning tree
and Steiner tree on random graphs.  There were also some earlier
investigations of the minimum bounded-diameter spanning tree on random
graphs \cite{MR1744176,MR1900907} which consisted
of testing the performance of several other heuristic algorithms.
In the present paper we rigorously
analyze the asymptotic weight of these trees. Indeed, we also describe
two algorithms that approximate well the constrained spanning (or
Steiner) tree problem.

\subsection{Proof strategy}

In the remaining sections of this paper, we present the proofs of
these theorems.  We saw already that \tref{MST-cutoff} is
implied by Theorems~\ref{Steiner-deep} and~\ref{th:k-small}.

The upper bound in \tref{th:k-small} follows from analyzing (in
\sref{greedy}) the Steiner tree produced by a simple
greedy heuristic algorithm, which estimates how many vertices should
occur in each level of the tree, and picks the cheapest set of this many
vertices connected to the previous level of the tree.  The proof of
the upper bound in \tref{Steiner-deep} is also algorithmic, and
appears in \sref{hybrid}.  There, the strategy is to
start with the minimum unconstrained spanning or Steiner tree, delete a
small number of edges to break the tree apart into pieces, and then
splice these pieces back together using the greedy algorithm from
\tref{th:k-small}.  The resulting tree has almost the same
weight as the original tree, and it has very small depth.  Both of
these algorithms produce trees rooted at a pre-specified vertex, so
they yield both the bounded-depth and bounded-diameter upper bounds.

The lower bounds in \tref{Steiner-deep} are self-evident, since
the weight of the unconstrained minimum Steiner tree is an obvious
lower bound on the weight of a constrained Steiner tree.  The lower
bounds in \tref{th:k-small} are proved in \sref{pf:lb},
and make use of a tight concentration inequality that is derived in
\sref{concentration}.  
The lower bound applies to any Steiner tree (where any vertex can be
the root in a tree of diameter $2k$, or any edge can be the root in a
tree of diameter $2k+1$), so recalling~\eqref{diam<depth}, the proof
yields both the bounded-diameter and bounded-depth lower bounds.

When carrying out the above calculations in \sref{greedy},
\sref{hybrid}, \sref{concentration}, and \sref{pf:lb}, we assume
that the edge weights are distributed according to exponential random
variables with mean~$1$, since this distribution is quite natural and
simplifies many of the calculations.  In \sref{distributions}, we show
how our results for exponential random variables imply the corresponding
results for other distributions (when the depth bound is bigger than~$1$).

\section{Greedy tree}\label{greedy}

\subsection{Construction} \label{construction}

Let us consider the following greedy method for algorithmically
growing a low-weight spanning or Steiner tree with bounded depth or
diameter.  We will build a tree in which every vertex is within
distance~$k$ from a particular root vertex or root edge.  (For
bounded-diameter trees, the root may be chosen arbitrarily.)  This
greedy tree~$\T_\ell$ depends on a sequence of nonnegative
numbers $\ell=(\ell_0,\ell_1,\dots,\ell_k)$ such that $\ell_0=1$ (if
the root is a vertex) or $\ell_0=2$ (if the root is an edge), and
$\sum_i\ell_i = n$ if the desired tree is a spanning tree, and
$\ell_k=m$ if the desired tree is a Steiner tree.  The idea is that
the greedy tree~$\T_\ell$ will have exactly $\ell_i$ vertices
at distance $i$ from the root, except that for Steiner trees there may
be fewer than $\ell_k=m$~vertices at the $k$\th level if the terminal
vertices were used closer to the root.  The construction is inductive.
Level~$0$ is the root vertex, or pair of vertices if the root is an
edge.  For convenience, we let $s_i=\ell_0+\dots+\ell_i$.  For $i < k$,
suppose we have chosen the $s_{i-1}$ vertices at level~$i-1$ and
below.  For each of the unchosen $n-s_{i-1}$ vertices, we look at the
lightest edge connecting it to level~$i-1$, and, to form level~$i$ of
the tree~$\T_\ell$, we choose the $\ell_i$ of these vertices
that have the lightest edges to level~$i-1$.  For level~$i=k$, we
connect the terminal vertices (that have not already been included
in the tree) using their lightest edge to level~$k-1$.

\subsection{Approximate weight}

The choice of the sequence $\ell$ has great influence on the
total weight of the resulting tree~$\T_\ell$.  Let $\wt_i(\T_\ell)$ be
the (random) total weight of the edges within the greedy tree
connecting level~$i-1$ to level~$i$.  For each of the $n-s_{i-1}$
vertices not in the tree up to level~$i-1$, the weight of the lowest
weight edge leading to it from the $\ell_{i-1}$ vertices at
level~$i-1$ is an exponential random variable with mean
$1/\ell_{i-1}$.  When picking the $\ell_i$ vertices at level~$i<k$,
we pick the $\ell_i$ smallest of these random variables. The $j$\th
smallest has expectation
\[
\frac{1}{\ell_{i-1}} \left[\frac{1}{n-s_{i-1}} + \frac{1}{n-s_{i-1}-1} +
  \dots + \frac{1}{n-s_{i-1}-(j-1)}\right],
\]
and so
\[
\E[\wt_i(\T_\ell)] = \frac{1}{\ell_{i-1}} \left[\frac{\ell_i}{n-s_{i-1}} +
  \frac{\ell_i-1}{n-s_{i-1}-1} + \dots +
  \frac{1}{n-s_{i-1}-(\ell_i-1)}\right].
\]
In \sref{concentration} we derive a concentration result for these
random variables.

In the case of spanning trees ($m=n$), the above formula also holds for level
$i=k$, and simplifies to $\E[\wt_k(\T_\ell)] = \ell_k/\ell_{k-1}$.
For general Steiner trees, at level~$k$ we have $\E[\wt_k(\T_\ell)]
\leq \ell_k/\ell_{k-1}$ because some of the terminals may have been
selected already.

Thus the expected weight $\E[\wt_i(\T_\ell)]$ of the $i$\th level
may be approximated by
\[
\E[\wt_i(\T_\ell)] \approx \begin{cases}
  \dfrac{\ell_i^2}{2 n \ell_{i-1}}, &   \text{if $s_i\ll n$;} \\
  \dfrac{\ell_i}{\ell_{i-1}} = \dfrac{\ell_i^2}{m \ell_{i-1}},
        & \text{if $s_i=n$,}\ \ \ \text{i.e., $i=k$.}
\end{cases}
\]
If $s_{k-1}\ll n$ then the above approximation
holds for all $i$.  Next we choose a good sequence~$\ell$
that makes $\E[\wt(\T_\ell)] = \sum_i \E[\wt_i(\T_\ell)]$ small.

\subsection{An optimization problem} \label{opt}

It is convenient to define $f_n(a,b) = b^2/(2 n a)$ and
\begin{align*}
 f_{n,c}(\ell) = f_{n,c}(\ell_0,\ell_1,\dots,\ell_k) & =
f_n(\ell_0,\ell_1)+\cdots+f_n(\ell_{k-2},\ell_{k-1})+c f_n(\ell_{k-1},\ell_k)\\
& =
\frac{\ell_1^2}{2 n \ell_0} + \cdots + \frac{\ell_{k-1}^2}{2 n \ell_{k-2}}
 + c \frac{\ell_k^2}{2 n \ell_{k-1}},
\end{align*}
(a factor of $c$ appears only in the last level).  We
have argued that the greedy tree with level sizes $\ell$ has expected
weight approximately $f_{n,2n/m}(\ell)$, provided that most of the nodes
occur in the last level.  Next we optimize $f_{n,c}(\ell)$, which is a
deterministic function of~$\ell$, allowing the level sizes to be real
numbers rather than constraining them to be integers.  We further
relax the constraints on the sum of the level sizes
($\ell_0+\cdots+\ell_k=n$ for MST and the corresponding constraint for
the Steiner tree), and instead fix $\ell_0(=1)$ and $\ell_k$, and
optimize the intermediate level sizes.
In \sref{small-depth} and \sref{small-diam} we return to the
question of how well this approximates the weight of the greedy tree
with constrained integer level sizes and random edge weights.

For $\ell_{i-1}$ and $\ell_{i+1}$  held fixed, let us find the choice of
$\ell_i$ which minimizes $f_n(\ell)$.
If $i+1<k$ then the two terms of $f_n(\ell)$ involving $\ell_i$ are
\[
f_n(\ell_{i-1},\ell_i) + f_n(\ell_i,\ell_{i+1}) =
\frac{\ell_i^2}{2 n \ell_{i-1}} +
\frac{\ell_{i+1}^2}{2 n \ell_{i}} ,
\]
which is minimized when
\[
\frac{2 \ell_i}{2 n \ell_{i-1}} - \frac{\ell_{i+1}^2}{2 n \ell_{i}^2} = 0,
\]
i.e., when $f_n(\ell_i,\ell_{i+1}) = 2 f_n(\ell_{i-1},\ell_i)$.
In the last level, with $i=k-1$ we still have that the optimal
choice of $\ell_i$ yields $f_n(\ell_i,\ell_{i+1}) = 2 f_n(\ell_{i-1},\ell_i)$.

An optimal sequence of $\ell_i$'s should satisfy this for all $i$, so we
wish to solve the recursion subject to $\ell_0=1$ with given $\ell_k$.
Let $r_i = \ell_i/\ell_{i-1}$.
If $i+1<k$ we have $ 2\ell_i^3 = \ell_{i-1} \ell_{i+1}^2$, i.e.,
\[
2 r_i = r_{i+1}^2.
\]
For $i=k-1$ we have $ 2\ell_i^3 = c \ell_{i-1}\ell_{i+1}^2$ so
\[
2 r_i = c r_{i+1}^2.
\]
When the ratios $r_i$ satisfy these equations, we have, for all $i<k$,
\[
r_i = 2(r_k \sqrt{c} / 2)^{2^{k-i}}.
\]
Multiplying, we find for $i<k$
\[
\ell_k/\ell_i = \prod_{j=i+1}^k r_j =
2^{k-i}(r_k \sqrt{c} / 2)^{2^{k-i}-1} / \sqrt{c}.
\]
In particular,
\[
\ell_k = \ell_k/\ell_0 =
2^{k} (r_k \sqrt{c}/2)^{2^{k}-1} / \sqrt{c},
\]
so
\[
r_k = \frac{2}{\sqrt{c}} \left(\frac{\ell_k \sqrt c}{2^{k}}\right)^{\frac1{2^k-1}}
\]
and for $i<k$,
\begin{equation}\label{li}
\ell_i = 2^{i} \left(\frac{\ell_k \sqrt c}{2^{k}}\right)^{1-\frac{2^{k-i}-1}{2^k-1}}.
\end{equation}
On this optimal sequence $\ell$, we have
\begin{equation}\label{fcl}
 f_{n,c}(\ell) = (2-2^{1-k}) \frac{c \ell_k^2}{2 n \ell_{k-1}}
 = \frac{2 \ell_k \sqrt c}{n} (1-2^{-k})
   \left(\frac{\ell_k \sqrt c}{2^{k}}\right)^{\frac1{2^k-1}}.
\end{equation}

\subsection{Bounded-depth trees} \label{small-depth}

How do we relate this sequence to the greedy minimum bounded-depth
spanning tree or Steiner tree?  Let $\hat\ell$ denote this optimal
sequence when $\ell_k = m$ and $c=2n/m$.  Let us greedily place
$\lceil \hat\ell_i \rceil$ nodes in level~$i$ for $i<k$, and then
connect the remaining nodes to the last level of the tree, and let
$\T_{\hat\ell}$ denote this tree.  For fixed
$m/n$, so long as $k \leq \log_2 \log n - \omega(1)$, we have $r_k \gg
1$, so for each $i<k$ we have $s_i\ll n$, so the expected weight of
level~$i$ of the tree is
$(1+o(1))f_n(\lceil\hat\ell_{i-1}\rceil,\lceil \hat\ell_i\rceil)$ for
$i<k$, and is at most $c f_n(\lceil\hat\ell_{i-1}\rceil,\lceil
\hat\ell_i\rceil)$ for $i=k$.  Furthermore, for each $i>0$ we have
$\hat\ell_i\gg 1$, so the rounding to integers only causes a
$(1+o(1))$ multiplicative correction.  Thus the expected weight of
this greedy minimum spanning tree or Steiner tree is at most
\begin{multline}
\E[\wt(\T_{\hat\ell})] \leq
(1+o(1)) f_{n,2n/m}(\hat \ell)
 = (1+o(1)) \sqrt{8m/n} (1-2^{-k})
            \left(\frac{\sqrt{2 m n}}{2^{k}}\right)^{\frac1{2^k-1}}\\
 \text{(for $k \leq \log_2 \log n - \omega(1)$)}.
\end{multline}
This is the upper bound of the in-expectation part of \tref{th:k-small}
for bounded-depth Steiner trees when the edge weights are $\Exp(1)$
random variables.  Concentration will follow when we prove the lower
bound in \sref{pf:lb}.

For example, when $k=2$ and $m=\alpha n$, the best choice is $\ell_1 \approx
\alpha^{1/3}n^{2/3}$, yielding a total expected weight of about
$\frac32\alpha^{2/3} n^{1/3}$.

We will also be interested in taking larger $k$'s, namely
$k=\log_2\log n -\Theta(1)$ and larger.  If we simply substitute this
$k$ (or any larger $k$) into the estimate for the weight of the greedy
tree, we would get $\E[\wt(\T_{\hat\ell})] \leq \Theta(\sqrt{m/n})$.
This estimate for $\E[\wt(\T_{\hat\ell})]$ is not valid, because for
$i=k-\Theta(1)$ we have $s_i=\Theta(n)$, so the $\E[W_i]$ are larger
than the above formula gives.  However, these $\E[W_i]$'s are only a
constant factor larger than predicted, so we still have
\begin{equation}
  \E[\wt(\T_{\hat\ell})]\leq\Theta(\sqrt{m/n})\ \ \ \text{(for $k\geq \log_2\log n -\Theta(1)$)}.
 \label{greedy-deep}
\end{equation}

Intuitively, all the choices that we made when constructing the greedy
tree are close to optimal.  We will see in the next section that it is
in fact possible to build a better tree by making non-greedy choices
when $k$ is larger than $\log_2\log n + \omega(1)$ (the
construction there makes use of this greedy tree, combining it with
the optimal minimum spanning tree with unbounded depth).  However, we
will prove in the lower bound section that the weight of this greedy
tree is within a factor of $(1+o(1))$ of the weight of the optimal
tree so long as $k\leq \log_2\log n - \omega(1)$.

\subsection{Bounded-diameter trees} \label{small-diam}

Having proved the upper bound for $\wt\!\Big(\MSTr{k}(K_n,m)\Big)$, we now
consider the case of trees with a bounded diameter. As noted, the primary
difference is that now there is no fixed root from which to measure
distances. The following is an easy observation (which will be more useful
for the lower bounds).

\begin{lemma}
  A tree with diameter $2k$ contains a unique root vertex, from which
  the tree has depth~$k$. A tree of diameter $2k+1$ contains a unique edge
  so that all vertices are within distance $k$ of an endpoint of the edge.
\end{lemma}

\begin{proof}
  Take a path of maximal length in the tree, and take as the root vertex or
  edge the central vertex or edge of the path. The bounds on the depth
  follow from the maximality of the path.
\end{proof}

To prove the upper bound for $\E\left[\wt\!\Big(\MSTd{2k}(K_n,m)\Big)\right]$,
just note that any tree of depth~$k$ is also a tree of diameter $2k$.  Thus
\[
\wt\!\Big(\MSTd{2k}(K_n,m)\Big) \leq \wt\!\Big(\MSTr{k}(K_n,m)\Big).
\]

To prove the upper bound for $\E\left[\wt\!\Big(\MSTd{2k+1}(K_n,m)\Big)\right]$, we fix
a root edge with small weight, then repeat the argument for the weight of
the greedy tree, except that level~$0$ now has size $\ell_0=2$ rather than 1.
There is an easy way to relate the optimal costs without repeating the
optimization problem.

The key is that the cost $f_{n,c}(\ell)$ is homogeneous in the sequence
$\ell$.  Consider the optimal sequence $\hat\ell$ for spanning $m/2$
terminals among $n/2$ vertices, then $2\hat\ell$ is the optimal
sequence for spanning $m$ terminals among $n$ vertices, except that it
starts with two vertices at level~$0$.  Thus we can repeat the greedy
construction to find that
\begin{align*}
\E\left[\wt\!\Big(\MSTd{2k+1}(K_n,m)\Big)\right] &\leq
(1+o(1)) f_{n,2n/m}(2\hat \ell) \\
 &= (1+o(1)) f_{n/2,2n/m}(\hat \ell) \\
  &= \text{our upper bound on $\E\left[\wt\!\Big(\MSTr{k}(K_{n/2},{m/2})\Big)\right]$.}
\end{align*}
This is the upper bound of the in-expectation part of \tref{th:k-small}
for bounded-diameter Steiner trees when the edge weights are $\Exp(1)$
random variables.

\section{Sliced-and-spliced tree}\label{hybrid}

We now show how to construct a small-diameter spanning (or Steiner)
tree with weight close to $\zeta(3)$ (or the weight of the
unconstrained Steiner tree).  The idea is to take the true
(unconstrained) minimum spanning tree (or Steiner tree), break it
apart into small subtrees which each still contain many vertices, and
then splice the subtrees together using the greedy tree approach from
\sref{greedy}.  The resulting tree is locally similar to the minimum
spanning / Steiner tree, so the weight is about the same, but globally
it has been rewired to have much smaller diameter.
\vspace*{3pt}

For the slicing part we use the following lemma:
\begin{lemma}\label{lem:slice}
  Any tree with diameter at least $\lfloor\Delta/2\rfloor$ contains a
  forest on the same vertex set where the components have diameters
  between $\lfloor\Delta/2\rfloor$ and $\Delta$.
\end{lemma}

Note that the weight of the resulting forest is bounded by the weight of
the tree.

\begin{proof}
  Any tree with diameter greater than $\Delta$ may be broken up into
  two trees with diameter at least $\lfloor \Delta/2\rfloor$ by removing
  the middle edge of some path realizing the diameter.  This may be
  iterated as long as there are components with diameter greater than
  $\Delta$.
\end{proof}

For the splicing part we use the following lemma:
\begin{lemma}\label{lem:partition}
  Suppose we are given a partition of the vertices of the complete
  graph $K_n$ into clusters of size at least $s$, and a root
  vertex~$r$.  Suppose the edges of $K_n$ are given independent
  $\Exp(1)$ edge weights.  If $\T$ is the minimal weight tree rooted
  at $r$ of depth at most $\log_2\log n$ that intersects each cluster
  of the partition, then $\E[\wt(\T)] = O(1/s)$.
\end{lemma}

\begin{proof}
  We repeat the greedy-tree construction of \sref{greedy} with some
  minor modifications.  We select the optimal level sizes
  $\ell_0,\dots,\ell_k$ (with $k\leq\lfloor\log\log n\rfloor$) for a
  spanning tree whose size is the number of clusters of the partition.
  During the construction, we will select a representative vertex from
  each cluster, and the constructed tree will contain only these
  representative vertices.  At step $0$ the partially constructed tree
  is the given root vertex.  At the $i$\th step, each cluster which is
  already connected to the partially constructed tree is connected
  through its representative vertex, while other clusters do not yet
  have a representative chosen.  For each unconnected cluster, the
  lightest edge from it to a representative vertex in level $i-1$ has
  weight which is dominated by $\Exp(1/s)$, and we choose the cheapest
  $\ell_i$ unconnected clusters to connect to level $i-1$.  Since the
  greedy tree of depth $\log_2\log (\text{number of clusters})$ has
  expected weight $O(1)$, it follows that the constructed tree~$\T$
  has expected weight $\E[\wt(\T)] = O(1/s)$.
\end{proof}

Recall that the weight of the unconstrained Steiner tree is
$(1-o(1))(m-1)/n \log(n/m)$ (w.h.p.\ and in expectation) when $m\ll n$
\cite{BGRS}, and that (using also Frieze's result on spanning trees
\cite{Frieze1985}) consequently $\wt(\MST(K_n,m))= \Theta(m/n \log(e n/m))$
for $2\leq m\leq n$.

\begin{theorem}[sliced-and-spliced tree]\label{deep-decay}
  Suppose $2\leq m \leq n$ and $k = \log_2\log n + \Delta$ where $\Delta
  \geq n/(m\log(e n/m))$, and the edge weights are exponential random
  variables with mean $1$. Then
  \[
  \E\left[\wt\!\Big(\MSTr{k}(K_n,m)\Big) - \wt(\MST(K_n,m))\right] \leq
  O\left(\sqrt{\frac{m \log(e n/m)}{n \Delta}}\right).
  \]
\end{theorem}

The idea for the construction of the sliced-and-spliced tree is to
take the minimum spanning / Steiner tree, slice it apart according to
\lref{lem:slice} into pieces of diameter $\Theta(\Delta)$, and splice
the pieces together as in \lref{lem:partition} using edges of total
weight $O(1/\Delta)$.  This does not quite work, since
\lref{lem:partition} assumes that the edge weights do not depend on
the partition, but this independence issue can be overcome by starting
with a slightly different tree as follows.

\begin{proof}[Proof of~\tref{deep-decay}]
  Each edge weight $w_e$ has $\Exp(1)$ law. For some $0<\eps<1/2$
  (eventually we use $\eps\ll 1$), we can write $w_e = \min(w'_e,
  w''_e)$, where $w'_e$ and $w''_e$ are independent exponentials with
  mean $1/(1-\epsilon)$ and $1/\epsilon$ respectively.  Let $\T'$ be
  the Steiner tree for the weights $w'$, and note that the weights
  $w''$ are independent of $\T'$.  The weight $\wt'(\T') \eqd
  \wt(\MST(K_n,m))/(1-\eps)$, where $\wt'$ is the weight using
  $w'$.  This implies $\E[\wt'(\T')] = (1+O(\eps)) \E[\wt(\MST(K_n,m))]$.

  If $\T'$ has diameter at most $\Delta$, then this is the
  sliced-and-spliced tree.  Otherwise, by \lref{lem:slice} we can
  slice $\T'$ into a forest $\mathcal F'$ whose connected components
  have size at least $\Delta/2$ and diameter at most $\Delta$.  Next,
  using Lemma~\ref{lem:partition} with the edge weights $w''$, the
  minimal tree $\T''$ of depth at most $\log_2\log n$ that connects
  the trees of $\mathcal F'$ has expected weight $O(1/(\eps\Delta))$.
  The sliced-and-spliced tree is then $\widehat \T = \mathcal{F}' \cup \T''$.
  Since any vertex is at distance at most $\Delta$ from $\T''$, it
  follows that $\widehat\T$
  has depth at most $\log_2\log n + \Delta$.  Moreover,
  \begin{align*}
  \E[\wt(\widehat\T)] &\leq \E[\wt(\T')] + \E[\wt(\T'')] \\
  &= \E[\wt(\MST(K_n,m))] + O(\eps \E[\wt(\MST(K_n,m))]) + O(2/(\eps\Delta)).
  \end{align*}

  It remains to pick $\eps$ to minimize this bound, namely $\eps =
  \Theta(1/\sqrt{\Delta \E[\wt(\MST(K_n,m))]})$.  Since $\Delta\geq
  n/(m\log(e n/m)) = \Theta(1/\E[\wt(\MST(K_n,m))])$ was one of our
  assumptions, we can pick such an $\eps\leq 1/2$.  Thus,
  $\wt\!\Big(\MSTr{k}(K_n,m)\Big)$ is at most the weight of the
  sliced-and-spliced tree, which is at most $\wt(\MST(K_n,m))$ plus a
  quantity which in expectation is at most $\Theta\left(\sqrt{\frac{m
        \log(e n/m)}{n \Delta}}\right)$.
\end{proof}

The above proof is wasteful in the separation of weights into two
independent components.  If this is not done, then the edges between
sub-trees are likely not to have come from the MST, and so tend to be
heavier.  It is plausible that this only increases the weight of the
connecting tree $\T''$ by a constant factor rather than a factor
of $O(1/\eps)$.

\begin{proof}[Proof of \tref{Steiner-deep} for exponential weights.]
  Using \tref{deep-decay}, if $k=\log_2\log n + \Delta$ where
  $\Delta \geq \omega(n/(m \log (e n/m)))$, then
  $$ \E\left[\wt\!\Big(\MSTr{k}(K_n,m)\Big) - \wt(\MST(K_n,m))\right]
  \leq o(\E[\wt(\MST(K_n,m))]), $$
  so the convergence in probability for bounded-depth Steiner trees is
  an immediate consequence of Markov's inequality, and convergence in
  expectation is also immediate.  The bounded-diameter statements are
  a consequence of the bounded-depth statements.
\end{proof}
The case of other distributions is handled in \sref{distributions}.

\section{Concentration of level weights}\label{concentration}

Let $U_1,\dots,U_p$ be a \underline{p}ool (set) of $p$ i.i.d.\ 
exponential random variables with mean $1$, and for $b\leq p$, let $W_{b,p}$ be the
sum of the $b$ \underline{b}est (smallest) $U_i$'s.  The distribution
of $W_{b,p}$ plays a key role in the behavior of the weights of
bounded-depth minimum spanning trees and bounded-depth Steiner trees.
The total weight~$W_i$ of the edges connecting levels $i-1$ and $i$
in the greedy tree from \sref{construction} is given by
$$
\wt_i(\T_\ell) = \frac{1}{\ell_{i-1}} W_{\ell_i,n-\ell_0-\ell_1-\dots-\ell_{i-1}}.
$$
We derive here some basic properties of $W_{b,p}$, including its
expected value, and the probability that it deviates far from its
expected value.

Let $Y_i$ be the $i$\th smallest of the $U_i$'s.  Since the minimum of
independent exponentials is again an exponential, and since an
exponential conditioned to be larger than some value is a translated
exponential, we have $Y_{i+1}-Y_i = \frac1{p-i} X_i$, where the
$X_i$'s are i.i.d.\ $\Exp(1)$ random variables (where by convention
$Y_0=0$).  It follows that
  \[ W_{b,p} =
  \sum_{i=1}^{b} Y_i = \sum_{i=0}^{b-1} \frac{b-i}{p-i} X_i.
  \]
Thus
$$
\E[W_{b,p}] = \sum_{i=0}^{b-1} \frac{b-i}{p-i}.
$$
Let us approximate the expected value by
$$ \overline W_{b,p} = \int_0^b \frac{b-i}{p-i} \, di
 = b + (p-b) \log\left(1- \frac{b}{p}\right)
 = \frac{b^2}{2 p} + \frac{b^3}{6 p^2} + \frac{b^4}{12 p^3} + \cdots 
 \begin{cases}\leq b^2/p, \\ \geq b^2/(2p);\end{cases}
$$
we have
$$
\overline W_{b,p} \leq \E[W_{b,p}] \leq \frac{b}{p} + \overline W_{b,p}.
$$

\begin{lemma}\label{L:exp_ld}
  For any $\delta>0$ and $b\leq p$,
  \[
  \Pr\left[W_{b,p} < (1-\delta) \frac{b^2}{2p} \right]
  \leq \Pr\left[W_{b,p} < (1-\delta) \E[W_{b,p}] \right]
  \leq \exp\left[-\frac18 \delta^2 b \right].
  \]
\end{lemma}

\begin{proof}
 We use the method of bounded differences (see e.g., \cite{MR1036755}).
 For $\beta>0$, we have
  \begin{align*}
    \Pr[W_{b,p} < x]
    = \Pr[e^{-\beta W_{b,p}} > e^{-\beta x}]
    &\leq e^{\beta x} \E[e^{-\beta W_{b,p}}] \\
    &= e^{\beta x} \prod_{i=0}^{b-1} \E[e^{-[\beta (b-i)/(p-i)]X_i}] \\
    &= e^{\beta x} \prod_{i=0}^{b-1} \frac{1}{1+\beta (b-i)/(p-i)} \\
    &= \exp \left[\beta x - \sum_{i=0}^{b-1}
      \log \left(1+\beta \frac{b-i}{p-i}\right) \right]. \\ 
    \intertext{Because $-\log(1+u)\leq -u+u^2/2$ for $u>0$, we have}
    \Pr[W_{b,p} < x]
    &\leq \exp \left[\beta x+\sum_{i=0}^{b-1}
      \left(-\beta \frac{b-i}{p-i} +
        \frac{\beta^2}{2} \frac{(b-i)^2}{(p-i)^2} \right) \right]\\
    &\leq \exp \left[\beta (x - \E[W_{b,p}]) + \frac{\beta^2}{2}
      \frac{b^3}{p^2} \right].
\intertext{Letting $x=(1-\delta) \E[W_{b,p}]$, we obtain}
    \Pr[W_{b,p} < (1-\delta) \E[W_{b,p}]]
  &\leq \exp \left[ - \beta \delta \E[W_{b,p}] +
    \frac{\beta^2}{2}\frac{b^3}{p^2} \right], 
\intertext{and setting $\beta = \delta \E[W_{b,p}] p^2 / b^3$, we obtain}
    \Pr[W_{b,p} < (1-\delta) \E[W_{b,p}]]
  &\leq \exp \left[ - \delta^2 \E[W_{b,p}]^2 \frac{p^2}{2b^3} \right]
\intertext{and since $\E[W_{b,p}] \geq b^2/(2p)$, we conclude}
    \Pr[W_{b,p} < (1-\delta) \E[W_{b,p}]]
  &\leq \exp \left[ - \delta^2 \frac{b}{8} \right]. \qedhere
  \end{align*}
\end{proof}

When bounding $\Pr[W_{b,p}<(1-\delta)b^2/(2p)]$, it is possible
to get a constant of $3/8$ rather than $1/8$, and $3/8$ is tight.  But
we also use $\Pr[W_{b,p}<(1-\delta)\E[W_{b,p}]]$ in
\sref{pf:lb}, and in the end this constant does not affect the
asymptotic lower bound that we prove there.

\section{MST lower bounds}\label{pf:lb}

\subsection{Strategy}

For sets~$A$ and $B$ of vertices, let $F(A,B)$ be the minimal total weight of a
set of edges connecting each vertex in $B$ to some vertex in $A$.
Note that $F(A,B)$ is increasing in $B$ and non-increasing in $A$. Define
\[
F(a,b) =
\min_{\substack{ |A|\leq a, |B|\geq b\\A\cap B = \emptyset}}  F(A,B) =
\min_{\substack{ |A|  =  a, |B|  =  b\\A\cap B = \emptyset}}  F(A,B)
,
\]
i.e., the minimal cost for connecting at least $b$ vertices to at most
$a$ vertices.  Next let
$$ F(\ell) = F(\ell_0,\ell_1,\dots,\ell_k) = \sum_{i=1}^k F(\ell_{i-1},\ell_i);$$
this is a lower bound on the cost of any spanning tree of depth~$k$
whose level sizes are given by $(\ell_0,\ell_1,\dots,\ell_k)$.
(Mnemonically, the $F$'s are random variables determined by the edge
weights of the graph, and the $f$'s from \sref{opt} are deterministic
quantities which we argue are likely to closely approximate the
$F$'s.)  We can obtain a sharper lower bound by treating the last
level differently.  In particular, for the last level of Steiner
trees, we need only consider sets~$B$ which contain only terminal
nodes of the Steiner tree.  The sharper bound is then
$$ F_m(\ell) = \sum_{i=1}^{k-1} F(\ell_{i-1},\ell_i) + F_m(\ell_{k-1},\ell_k),$$
where $F_m(a,b)$ is defined as $F(a,b)$ was, but with the set $B$
restricted to be a subset of the $m$ terminals of the Steiner tree.

Let $(\hat\ell_0,\hat\ell_1,\dots,\hat\ell_k)$ be the ``greedy
sequence'' of level sizes that optimizes $f_{n,2n/m}(\hat\ell)$ and which we
used for the greedy tree in \sref{greedy}.  Our strategy is to show that
$F_m(\ell)$ is approximately minimized at $F_m(\hat\ell)$ for $n$
large enough, and that $F_m(\hat\ell)$ is within a factor $(1-\delta)$
of $f_{n,2n/m}(\hat\ell)$.  It will follow that for $n$ large enough, the
weight of the greedy tree is close to the weight of the optimal tree.

\begin{lemma}\label{L:no_cheap_sets}
  For any $\delta > 0$, $a \in \{1,\ldots,n\}$ and $b \in \{1,\ldots,
  n-a\}$ we have
  \[
  \Pr\left[ F(a,b) < (1-\delta) \frac{b^2}{2 n a} \right]
  \leq \exp\left[a\log\frac{n e}{a} - \frac18 \delta^2 b\right].
  \]
\end{lemma}

\begin{proof}
  Fix a set $A$ of size $a$.  For each $x\notin A$ the minimal
  weight of an edge connecting $x$ to $A$ is an independent $\frac 1{a}
  \Exp(1)$.  Let $W_A$ be the total weight of edges connecting the cheapest $b$
   vertices to $A$; $W_A$ and $\frac{1}{a} W_{b,n-a}$ have the same
   distribution.  By \lemref{exp_ld},
  \[
  \Pr\left[W_A < (1-\delta)\frac{b^2}{2na}\right] \leq
  \Pr\left[ W_A < (1-\delta) \frac{b^2}{2(n-|A|)|A|} \right]
  \leq \exp\left[-\frac18 \delta^2 b \right].
  \]
  Finally, the number of sets $A$ of size $a$ is
  $$
  \binom{n}{a}
  \leq \left(\frac{n e}{a}\right)^a, $$
  and a union bound yields the claim.
\end{proof}

\begin{lemma}\label{L:no_cheap_sets2}
  For any $\delta > 0$, $a \in \{1,\ldots,n\}$ and $b \in \{1,\ldots,
  m-a\}$ we have
  \[
  \Pr\left[ F_m(a,b) < (1-\delta) \underbrace{\left(1-\frac{m-b}{b}\log
        \frac{m}{m-b}\right)}_{\text{$1-o(1)$ if $m/b\to1$}} \frac{b}{a}
  \right] 
  \leq \exp\left[a\log\frac{n e}{a} - \frac18 \delta^2 b\right].
  \]
Note that the bound on $F_m(a,b)$ is $(1-o(1)) b^2/(m a)$
in the limit $\delta\to0$ and $b/m\to1$.
\end{lemma}

\begin{proof}
  The proof is essentially the same proof used for
  \lref{L:no_cheap_sets}, except that the cheapest $b$ vertices
  come from a set of $m$ vertices, we use the second inequality
  from \lref{L:exp_ld} rather than the first and second combined,
  and we use the bound
  $$\E[W_A] = \frac{1}{a} \E[W_{b,m}]\geq \frac{1}{a} \overline W_{b,m} = 
\frac{b}{a} + \frac{m-b}{a} \log\left(1- \frac{b}{m}\right)
 = \left(1 - \frac{m-b}{b} \log \frac{m}{m-b}\right)\frac{b}{a}. \qedhere$$
\end{proof}

In other words, if $b$ is large enough compared to $a$, then it is unlikely
that $F(a,b)$ is much smaller than $f_n(a,b)=b^2/(2 n a)$, and $F_m(a,b)$
is unlikely to be much smaller than $f_{n,2n/m}(a,b) = b^2/(m a)$.
(Recall the definition of $f_{n,c}$ from \sref{opt}.)
Let us define
$$ \Rdn(a) = \frac{32}{\delta^2} a \log \frac{n e}{a}. $$
Let $$\fR_n(a,b) = \1_{b > \Rdn(a)} \frac{b^2}{2 n a}.$$
The reason for introducing this cutoff $\Rdn(a)$ is so that w.h.p.\
$F(a,b) \geq (1-\delta)\fR_n(a,b)$ regardless of what $b$ is; for a
given $a$ and $b$ the probability that this fails is at most
$\exp(-3 a \log(n e/a))$.  This bound is decreasing for $a$ in the
range $0\leq a \leq n e$, and since $1\leq a \leq n$, for any given
$a$ and $b$ the failure probability is at most $1/n^3$.  Upon summing
over the choices of $a$ and $b$, it follows that with probability at
least $1-1/n$, for all $a$ and $b$ we have $F(a,b) \geq
(1-\delta)\fR_n(a,b)$.

Let
$$ \fR_n(\ell) = \fR_n(\ell_0,\ell_1,\dots,\ell_k)
 = \sum_{i=1}^k \fR_n(\ell_{i-1},\ell_i)
 = \sum_{i=1}^k \frac{\ell_{i}^2}{2n\ell_i} \1_{\ell_i > \Rdn(\ell_{i-1})},$$
and $\fR_{n,c}(\ell)$ be defined similarly, but with an extra factor of $c$ in
the $k$\th term of the sum.

\begin{corollary}\label{C:no_cheap_2}
  With high probability ($\geq 1-1/n$), any spanning tree with level
  sizes given by $(\ell_0,\dots,\ell_k)$ has weight at least
  $(1-\delta) \fR_n(\ell)$.  If $\ell_k=(1-o(1)) m$, then any Steiner
tree connecting a given set of $m$ terminals with level sizes $(\ell_0,\dots,\ell_k)$ has weight at least
$(1-\delta-o(1)) \fR_{n,2n/m}(\ell)$.
\end{corollary}

Thus we are done if we show that $\fR_{n,c}(\ell)$ constrained to $\sum_i
\ell_i \geq m$ is almost minimized at the sequence $\hat \ell =
(\hat\ell_0,\dots,\hat\ell_k)$ that minimizes $f_{n,c}(\hat\ell)$
constrained to $\ell_k=m$ (and which we used in the greedy tree
construction).

\begin{definition} \label{large}
  We say that the $k$\th level of a sequence~$\ell$ is \textit{large} if
  $\ell_k \ge  (1-\delta^2/32)m$.  For $t\geq 1$, we say that the
  $(k-t)$\th level is \textit{large} if
  \[
  \ell_{k-t} \geq \overbrace{\Rdn^{-1}(\cdots \Rdn^{-1}}^t((1-\delta^2/32)m)\cdots).
  \]
\end{definition}

Since $\Rdn$ is monotone increasing in $a$ up to $a=n$, the inverse
function $\Rdn^{-1}(b)$ is well-defined for $1\leq b\leq n$.  Since
$\Rdn^{-1}(b)\leq (\delta^2/32)b$, and $\sum_i \ell_i \geq m$, it follows
that there must be at least one large level.

We may enlarge the set of $\ell$'s over which we are
minimizing; so long as $\fR_{n,c}(\ell)$ is still almost minimized at
$\hat\ell$, we will have our desired lower bound.  Naturally we relax
the constraint $\ell_i\in\N$ to $\ell_i\in\R^+$.  We keep the
constraint $\ell_0=1$.  We shall drop the $\sum_i
\ell_i \geq m$ constraint, and replace it with a constraint that there
is a large level in the above sense, since this only increases the set
of sequences that we are optimizing over.

\subsection{No small jumps}

Call a jump from $\ell_i$ to $\ell_{i+1}$ large if $\ell_{i+1}>
\Rdn(\ell_i)$ and small otherwise.  Small jumps contribute $0$ to
$\fR_{n,c}(\ell)$.  Suppose that a sequence contains a small jump
$(\ell_i,\ell_{i+1})$.  If $\ell_{i+1}<\Rdn(\ell_i)$, then we may
increase $\ell_{i+1}$ or decrease $\ell_i$, and each term of
$\fR_{n,c}(\ell)$ either stays the same or decreases.  Thus any small
jump in a sequence~$\ell$ minimizing $\fR_{n,c}$ is a jump from $a$ to
$\Rdn(a)$.  If two consecutive jumps of a sequence~$\ell$ minimizing
$\fR_{n,c}$ are large, then the intermediate value must satisfy
$2\ell_i^3=\ell_{i-1}\ell_{i+1}^2$.  We wish to show that a sequence
achieving the minimum value of $\fR_{n,c}$ in fact has no small jumps,
which will allow us to find the best sequence.  We start by showing
that it does not have a small jump followed by a large jump.

\begin{lemma} \label{SL}
  There is a $\delta_0>0$ so that whenever $c\geq 1$ and
  $0<\delta\leq\delta_0$, and a sequence $\ell$ has a small jump
  $(\ell_{i-1},\ell_i)$ followed by a large jump
  $(\ell_i,\ell_{i+1})$, it is possible to change $\ell_i$ so as to
  reduce $\fR_{n,c}(\ell)$.
\end{lemma}

\begin{proof}
  Let $\ell_{i-1}=a$ and $\ell_{i+1}=b$.  Let $C=32/\delta^2$.  Since
  $(\ell_{i-1},\ell_i)$ is a small jump, $\ell_i \leq
  \Rdn(\ell_{i-1})$.  If $\ell_i < \Rdn(\ell_{i-1})$, then we may
  replace $\ell_i$ with $\tilde\ell_i=\Rdn(\ell_{i-1})$ to get a new
  sequence $\tilde\ell$ for which
  $\fR_{n}(\tilde\ell_{i-1},\tilde\ell_i)$ is still $0$ but
  $\fR_{n}(\tilde\ell_{i},\tilde\ell_{i+1})<\fR_{n}(\ell_{i},\ell_{i+1})$,
  so $\fR_{n,c}(\tilde\ell) < \fR_{n,c}(\ell)$.  Next we consider the
  case $\ell_i = \Rdn(\ell_{i-1})$, and show that the sequence $\ell$
  can still be improved.  There is a slight difference when $i+1=k$,
  as opposed to $i+1<k$, since the last ($k$\th) summand of
  $\fR_{n,c}(\ell)$ contains a factor of $c$.  We deal below with the
  case $i+1=k$.  The case $i+1<k$ differs only in that $c$ does not
  appear, and is derived by replacing all $c$'s by 1's.  We consider
  the original sequence and two possible replacements of $\ell_i$ by
  $\ell'_i$ and $\ell''_i$ defined by
  \begin{align*}
    2(\ell'_i)^3 &= c \ell_{i-1}\ell_{i+1}^2, &
    \Rdn(\ell''_i) &= \ell_{i+1}.
  \end{align*}
  (The replacement $\ell'_i$ is optimal for two large jumps, and $\ell''_i$
  is optimal for a large jump followed by a small jump.)
  Let $U$, $U'$, and $U''$ be the contributions to $\fR_{n,c}$ from
  the two jumps in the three cases $(\ell_{i-1},\ell_i,\ell_{i+1})$,
  $(\ell_{i-1},\ell'_i,\ell_{i+1})$, and
  $(\ell_{i-1},\ell''_i,\ell_{i+1})$.  For $U$ and $U''$, by definition
  the jump that is small contributes $0$, so there is only one term.
  \begin{align*}
    U &= 0+c\frac{\ell_{i+1}^2}{2n\ell_i}
    = \frac{c b^2}{2n\Rdn(a)} = \frac{c b^2}{2 n C a\log(ne/a)}, \\
    U' &= \frac{{\ell'_i}^2}{2n\ell_{i-1}} + c \frac{\ell_{i+1}^2}{2n\ell'_i}
    = \frac{3}{2n} \left(\frac{c^2 b^4}{4a}\right)^{1/3},\\
    U'' &= \frac{(\ell''_{i})^2}{2n\ell_{i-1}}+0.
\intertext{
Since
$b = \Rdn(\ell''_i) = C \ell''_i \log (n e /\ell''_i) \geq \ell''_i$,
we have $b\geq C \ell''_i \log(n e/b)$, so
}
  U'' & \leq \frac{b^2}{2 n a C^2\log(ne/b)^2}.
\end{align*}
If $U''>U$ then
  \[
  C c \log(ne/b)^2 < \log(ne/a).
  \]
 If $U'>U$, then
  \[
  b < C' c^{-1/2} a \log(ne/a)^{3/2},
  \]
  where $C' = (3 C)^{3/2}/2 = 192\sqrt{6}/\delta^3$.

  If both $U'>U$ and $U''>U$ then we find
  \[
  \log(ne/a) > C c \log^2\frac{n e}{C' c^{-1/2} a\log(ne/a)}.
  \]
  If we denote $q=\log(ne/a)$, this becomes
  \[
  q > C c (q - \log(C' c^{-1/2} q))^2.
  \]
  
  Since $b\leq m$ and each jump increases the level size by at least
  $32/\delta^2$, we have $a\leq (\delta^4/1024)m$, so $q\geq \log[(1024 e /
  \delta^4) n/m]$.  Since $C'=192\sqrt6/\delta^3$, we have $q - \log C'
  \geq \frac14 \log q - \text{const}$.  Thus
  \[
  q > \frac{32}{\delta^2} c \left(\frac14 q - \log q +\frac12 \log c -
    \text{const}\right)^2. 
  \]
  But $q\to\infty$ as $\delta\to0$ and $c\geq 1$,
  so this equation cannot be true
  for small enough~$\delta$.  Thus, provided $\delta<\delta_0$, we
  have $\min\{U',U''\}\le U$, so replacing $\ell_i$ with one of
  $\ell_i'$ or $\ell_i''$ reduces $\fR_{n,c}(\ell)$.
\end{proof}

Consider the maximal $t$ such that level~$k-t$ of the sequence $\ell$ is
large. The jump from $\ell_{k-t-1}$ to $\ell_{k-t}$ must be a large jump,
since otherwise level~$k-t-1$ would also be a large level. Because (for
small enough $\delta$) there are no small jumps followed by large jumps, it
follows that the subsequence $\ell_0,\ell_1,\dots,\ell_{k-t}$ consists only
of large jumps.  If $t\neq0$ we need to show that this would imply that
the total cost up to level~$k-t$ is too high for $\ell$ to be
optimal.  Our next step is to obtain a lower bound on $\ell_{k-t}$, for
which we need the following lemma.

\begin{lemma}
Assume $\delta\leq 1$.  Recall that $C=32/\delta^2$, and that $\Gamma$ is the
gamma function (generalized factorial).
If for some $r\geq 2$
$$ b \geq \frac{e n}{(C \log C)^{r-1} \Gamma(r)^2},$$
then $$ \Rdn^{-1}(b) \geq \frac{e n}{(C \log C)^{r} \Gamma(r+1)^2} .$$
\end{lemma}

\begin{proof}
If $$a\leq \frac{e n}{(C\log C)^r \Gamma(r+1)^2},$$ then
\begin{align*}
\Rdn(a) &\leq C \frac{e n}{(C\log C)^r \Gamma(r+1)^2} \log \frac{n e}{e
  n/((C\log C)^r \Gamma(r+1)^2)} \\ 
 &\leq \frac{r\log(r^2 C\log C)}{r^2\log C} \times
       \frac{e n}{(C\log C)^{r-1}\Gamma(r)^2}.
\end{align*}
The first term is
\begin{align*}
\frac{r\log(r^2 C\log C)}{r^2\log C}
 &= \frac{2\log r}{r\log C} + \frac{\log C}{r \log C}
     + \frac{\log \log C}{r \log C}, \\
\intertext{and since $(\log r)/r\leq 1/e$ and $(\log\log C)/\log C\leq 1/e$
 and $r\geq 2$, we have}
 &\leq \frac{2}{e \log C} + \frac{1}{2} + \frac{1}{2 e}.
\end{align*}
As long as $\delta\leq 1$, we have $C\geq 32$, so that this quantity
is bounded by $1$.
\end{proof}

Thus we get a lower bound on the size $\ell_{k-t}$ of the first large level.

\begin{lemma} \label{lkt}
  If $\delta\leq1$ and $t\geq 1$ and level~$k-t$ is a large level, then
  $$ \ell_{k-t} \geq \frac{m}{(t \log\frac{n}{m} / \delta)^{\Theta(t)}}.$$
\end{lemma}
\begin{proof}
By our definition of ``large,''
  $$ \ell_{k-t} \geq \overbrace{R_{\delta,n}^{-1}(\cdots
    R_{\delta,n}^{-1}}^t((1-\delta^2/32)m)\cdots).$$ 
  Next we find the smallest~$r\geq 2$ satisfying
\begin{equation}\label{mnr}
  (1-\delta^2/32) m \geq \frac{e n}{(C \log C)^{r} \Gamma(r+1)^2},
\end{equation}
  the relevant $r$ satisfies $r\leq O(\log \frac{n}{m} / \log\log \frac{n}{m})$,
  so we can bound
\begin{align*}
  \ell_{k-t}& \geq \frac{e n}{(C \log C)^{r+t} \Gamma(r+t+1)^2} \\
\intertext{which, if $r>2$ (so that \eqref{mnr} is tight), can be bounded by}
 \ell_{k-t}&\geq  \frac{(1-\delta^2/32) m}{(C \log C(r+t)^2)^{t}}.
\end{align*}
Whether or not $r>2$, we can bound
\begin{align*}
\ell_{k-t} &\geq \frac{m}{(t \log\frac{n}{m}/\delta)^{\Theta(t)}}. \qedhere
\end{align*}
\end{proof}

\begin{lemma}\label{L:optimal_no_small_jumps}
  If $0<\delta<\delta_0$, there is a constant~$\Delta$ such that
  whenever $k \leq \log_2\log m - \log_2\log(e n/m) - \Delta$,
  any sequence~$\ell$ optimizing $\fR_{n,2n/m}(\ell)$ contains no
  small jumps.
\end{lemma}

\begin{proof}
If $t\geq 1$, then we may bound
$\fR_{n,c}(\ell)$ from below by the cost of the first $k-t$ levels, which by our
earlier calculation~\eqref{fcl} is
$$ f_n(\ell_0,\dots,\ell_{k-t}) = \frac{2 \ell_{k-t}}{n} (1-2^{t-k})
\left(\frac{\ell_{k-t}}{2^{k-t}}\right)^{\frac1{2^{k-t}-1}}. $$ 

Let us assume $t<k$. Upon substituting our lower bound for $\ell_{k-t}$
from \lref{lkt} and $c=2n/m$, we may compare this (LHS) to the
sequence from the greedy construction (RHS).  Since the expressions
are complicated, we do a sequence of inequality-preserving transformations
to determine which one is bigger:
\newcommand{\compare}{\stackrel{?}{\lessgtr}}
\begin{align*}
  \frac{2 m}{n (\frac{t}{\delta} \log\frac{n}{m})^{\Theta(t)}} (1\!-\!2^{t-k})
  \left(\!\frac{m/(\frac{t}{\delta} \log\frac{n}{m} )^{\Theta(t)}}
    {2^{k-t}}\!\right)^{\frac1{2^{k-t}-1}} &\compare
  \sqrt{\frac{8 m}{n}} (1-2^{-k})
  \left(\frac{\sqrt{2mn}}{2^{k}}\right)^{\frac1{2^{k}-1}} \\ 
 \left(\frac{m}{2^{k-t}}\right)^{\frac1{2^{k-t}-1}}
&\compare
\left(\frac{t}{\delta} \log\frac{n}{ m}\right)^{\Theta(t)}
\sqrt{\frac{n}{m}}
\left(\frac{\sqrt{2mn}}{2^{k}}\right)^{\frac1{2^{k}-1}} \\
 \left(\frac{m}{2^{k}}\right)^{\frac1{2^{k-t}-1}-\frac1{2^{k}-1}}
&\compare
\left(\frac{t}{\delta} \log\frac{n}{ m}\right)^{\Theta(t)}
\left(\sqrt{\frac{n}{m}}\right)^{\frac{2^k}{2^{k}-1}}\\
\Theta\!\left(\frac{2^t}{2^k}\right)(\log m - O(k))
&\compare
\Theta\!\Big(\!t\log \frac{t}\delta\Big) +
\Theta\!\Big(\!t \log\log\frac{n}{m}\Big) +
\Theta\!\left(\log\frac{n}{m}\right)\!.
\end{align*}
Let us assume $k\leq \log_2\log m - \Delta$, where $\Delta$ is a suitably
large constant depending on $\delta$. Then the $O(k)$ term in the LHS may
be neglected, and the LHS is $\Omega(2^\Delta 2^t)$, which is
larger than the $\Theta(t\log t)$ term on the RHS. If in addition, $k \leq
\log_2 \log m - \log_2 \log(e n/m) - \Delta$, then half the LHS is also
$\Omega(2^\Delta 2^t\log(n/m))$, which dominates the second and
third terms in the RHS.

Thus under these conditions on $k$, any sequence $\ell$ with $t\neq0$ is
not optimal, and so any optimizing sequence for $\fR_{n,2n/m}(\ell)$ contains
no small jumps.
\end{proof}

\subsection{Lower bounds}

\begin{proof}[Proof of \tref{th:k-small}, lower bounds, for exponential weights]
  We start with the bound on $\wt\!\Big(\MSTd{2k}(K_n,m)\Big)$. Combining
  \cref{C:no_cheap_2} with \lemref{optimal_no_small_jumps}, we find
  that w.h.p.\ the cost of a Steiner tree with level sizes given by $\ell$
  is at least $(1-\delta)\fR_{n,2n/m}(\ell)$, and that this is minimized by
  the sequence constructed in \sref{greedy} to give the upper
  bound. Thus we find the upper bound is tight.

  As noted, $\wt\!\Big(\MSTr{k}(K_n,m)\Big) \ge
  \wt\!\Big(\MSTd{2k}(K_n,m)\Big)$, so its bound is also tight.

  Finally substituting $m/2$ and $n/2$ throughout, we find that no tree can
  improve by more than $(1+o(1))$ on the greedy tree construction for the
  odd diameter case.
\end{proof}

\begin{proof}[Proof of \tref{th:k-small}, concentration, for exponential weights]
  For each of the random variables $\wt\!\Big(\MSTd{2k}(K_n,m)\Big)$,
  $\wt\!\Big(\MSTd{2k+1}(K_n,m)\Big)$, and
  $\wt\!\Big(\MSTr{k}(K_n,m)\Big)$, the random variable is almost
  never smaller than a factor of $1+o(1)$ smaller than our upper bound
  on their expected values.  It follows that for each of these random
  variables, the expected value is within a factor of $1+o(1)$ of our
  upper bound on it, and that these random variables are with high
  probability within a factor of $1+o(1)$ of their expected values.
\end{proof}

\section{Other weight distributions} \label{distributions}
\newcommand{\rup}[1]{#1^\eps}
\newcommand{\rdn}[1]{#1_\eps}
\newcommand{\apx}{\widetilde}

In the proofs up to this point, we assumed that the edge weights of
$K_n$ are distributed according to an exponential random variable with
mean~$1$.  In this section we prove the parts of
Theorems~\ref{MST-cutoff}, \ref{Steiner-deep}, and \ref{th:k-small}
that pertain to more general weight distributions.  In our notation up
until now we suppressed the weight distribution, but here we make it
more explicit: we let $K_n^{\apx W}$ denote the complete graph where
each edge weight is an i.i.d.\ copy of a non-negative random
variable~$\apx W$.

The key observation (which was made earlier in the context of
unconstrained minimum spanning trees \cite{Frieze1985,MR905183}
and Steiner trees \cite{BGRS}) is that w.h.p.\ 
only edges with weights~$o(1)$ are ever used (except when the depth
is~$1$).  Thus, it is principally the density of the distribution
near~$0$ that is significant.  
In this section we assume that the edge weights are i.i.d., and are
distributed according to some non-negative random variable $\apx W$
that has density~$1$ near~$0$, i.e., for positive~$t$ near~$0$,
$$\Pr[\apx{W}<t] = t + o(t).$$
(If the density near~$0$ exists and is not~$1$, then linearity in the
weights gives a multiplicative constant in the theorems.)  We let $W$
denote an exponential random variable with mean~$1$,
$$ W \sim\Exp(1),$$
and for $\eps>0$ define
\[
\rdn{W} = \begin{cases} W & W \leq \eps \\ \eps & W>\eps, \end{cases}
\ \ \ \ \ \ \ \text{and} \ \ \ \ \ \ \ 
\rup{W} = \begin{cases} W & W \leq \eps \\ \infty & W>\eps. \end{cases}
\]
For a generally distributed non-negative weight~$\apx{W}$ with
density~$1$ at~$0$, for any $\delta>0$ there is an $\eps>0$ such that
\[
(1-\delta)\rdn{W} \prec \apx{W} \prec (1+\delta)\rup{W},
\]
i.e., $\apx W$ is stochastically sandwiched between
$(1-\delta)\rdn{W}$ and $(1+\delta)\rup{W}$.  Since $\wt(\MST)$ is
monotone in the edge weights, it follows that bounded-depth/diameter
tree weight distributions are also stochastically sandwiched.

Let us call an edge of
the weighted graph \textit{$\eps$-light\/} if its weight is at most~$\eps$,
and otherwise let us call it \textit{$\eps$-heavy}.  The idea is to show that
w.h.p.\ the greedy, sliced-and-spliced, and optimal trees use only
light edges (when $k\geq 2$), so that it makes little difference
whether the edge weights are distributed according to~$W$ or~$\apx W$.

As we shall see, the density-$1$-at-$0$ assumption is enough to get
convergence in probability, but some additional assumption to rule out
the possibility of very fat tails in the distribution is required to
get the upper bounds for the convergence in expectation of the tree
weights.  We shall assume $$\E[\apx{W}]<\infty$$
when deriving convergence in expectation.

\subsection{Upper bounds}

We start by proving that, in the greedy tree and sliced-and-spliced
tree, heavy edges are rare.

\begin{lemma} \label{greedy-heavy}
  If $2\leq k\leq \frac14 \log_2 n$, the expected number of
  $\eps$-heavy edges contained within the greedy Steiner tree from
  \sref{greedy} is at most $\exp[\Theta(\log n) - \Theta(\min(\eps,1) n)]
  + k n \exp[-\eps n^{1/8}]$. 
\end{lemma}

\begin{proof}
  The level sizes~$\lceil\hat\ell_i\rceil$ from \eqref{li} are monotone
  increasing in~$i$, and monotone decreasing in~$k$.  At level~$1$,
  since $k\geq 2$ we have
  $$ (n/4^k)^{1/4} < 2 \sqrt[4]{2 n m / 4^k}\leq \hat\ell_1 \leq 2 \sqrt[3]{n m /2} < 2 n^{2/3}.$$
  The number of light edges emanating from the root is a binomial
  with parameters $n-1$ and $p=1-e^{-\eps}\geq\eps-\eps^2/2$.  A standard
  large-deviation formula (see \cite[Eqn.~5.6]{MR1036755}) tells us
  that for any binomial random variable $D$,
  $$\Pr[D<\E[D]/2]\leq e^{-\E[D]/8}.$$
  Assuming $\eps\geq 5 n^{-1/3}$ and $n$ is large,
  so that $\E[D]/2 \geq 2 n^{2/3}+1$,
  we deduce that the expected number of heavy edges in the first level
  of the greedy tree is at most
  $$
\lceil 2 n^{2/3}\rceil e^{-(n-1)p/8} = \exp[\Theta(\log n) - \Theta(p n)].
$$ 
  (If $\eps<5 n^{-1/3}$ or $n$ is not large,
  the conclusions of the lemma are trivially true.)
  For any subsequent level
  of the tree, the number of heavy edges is at most the number of
  vertices not connected to it via a light edge, and since there are
  at least $(n/4^k)^{1/4}$ vertices in the previous level, the
  expected number such vertices not reachable by a light edge
  is at most $n (1-p)^{(n/4^k)^{1/4}} =
  n \exp[-\eps (n/4^k)^{1/4}] \leq n \exp[-\eps n^{1/8}]$.  Upon
  multiplying by~$k-1$ (since there are $k-1$ levels after the first)
  and adding the heavy edges from the first level, we obtain the
  desired bound.
\end{proof}

\begin{proof}[Proof of \tref{th:k-small}, upper bounds, other weight distributions]
  The upper bounds for convergence in probability are an immediate
  consequence of the fact that for fixed~$\eps$, w.h.p.\ there are not
  any heavy edges.  The upper bounds for convergence in expectation
  follow from the fact the expected number of heavy edges is~$o(1)$,
  and the fact that $\E[\apx W]$ is finite.
\end{proof}

The following lemma essentially appears in \cite{BGRS}.
\begin{lemma} \label{Steiner-heavy}
  For $n\geq 3$ and $\eps>0$, the expected number of $\eps$-heavy
  edges in $\MST(K_n,m)$ is at most $O(e^{-\eps n} n^4\log^2 n)$.
\end{lemma}
\begin{proof}
  Consider any edge of the Steiner tree $\MST(K_n,m)$ with weight
  greater than $\eps$.  If its endpoints are connected by a path with
  total weight $\leq\eps$, then we could delete the heavy edge and
  replace it with a portion of or all of the low-weight path
  connecting that edge's endpoints, obtaining a lighter Steiner tree.
  (This argument appeared in \cite{BGRS}.)
  Janson proved that in the complete graph~$K_n$ with exponential edge
  weights, w.h.p.\ every pair of vertices is connected by a path of
  weight at most $(3+o(1)) n^{-1} \log n$ \cite{MR1723648}.  In fact,
  it follows from \cite[Eqn.~2.8]{MR1723648} that, when $n\geq 3$ and
  $\eps\geq 0$, the expected number of pairs of vertices not connected
  by a path of weight $\leq \eps$ is at most $O(e^{-\eps n} n^4\log^2
  n)$.  Thus, the expected number of heavy edges in $\MST(K_n,m)$ is
  at most $O(e^{-\eps n} n^4\log^2 n)$.
\end{proof} 

\enlargethispage{12bp}
\begin{lemma} \label{slice-n-splice-heavy}
  If $k=\log_2\log n + \Delta$ where $\Delta\geq n/(m\log(e
  n/m))$, the expected number of $\eps$-heavy edges in the
  sliced-and-spliced Steiner tree from \sref{hybrid} is at most
  $\exp[\Theta(\log n)-\eps n] + O(1/(\eps\sqrt\Delta))$.
\end{lemma}
\begin{proof}
  There could be as many as $\exp[\Theta(\log n)-\eps n]$ heavy edges
  in the starting Steiner tree $\MST(K_n,m)$.
  Of course, when we do the slicing of $\MST(K_n,m)$, no heavy edges
  are introduced, but heavy edges could be introduced when we splice
  the subtrees using the greedy-tree construction.  In the
  construction, recall that the total weight of the splice edges was
  in expectation at most 
  $$O\left(\sqrt{\frac{m \log(e n/m)}{n
        \Delta}}\right) \leq O(1/\sqrt{\Delta}).$$
  The expected number of $\eps$-heavy splice edges can be
  at most $1/\eps$ times as large as this.
\end{proof}

\begin{proof}[Proof of \tref{Steiner-deep}, other weight distributions]
  As above, the upper bounds for convergence in probability are an
  immediate consequence of the fact that for fixed~$\eps$, w.h.p.\ 
  there are not any heavy edges, and the upper bounds for convergence
  in expectation follow from the fact the expected number of heavy
  edges is~$o(1)$, and the fact that $\E[\apx W]$ is finite.
  
  The lower bounds follow from the fact that the unrestricted Steiner
  tree $\MST(K_n,m)$ w.h.p.\ has no heavy edges, and is at most as
  heavy as the bounded-depth/diameter Steiner trees.
\end{proof}

\subsection{Lower bounds}

\begin{proof}[Proof of \tref{th:k-small}, lower bounds, other weight distributions]
  Fix some $\delta>0$.
  Let $\rdn F(a,b)$ be defined as $F(a,b)$ but using $\rdn{ (W_e)}$.
  We argue that
  w.h.p., for every pair $a\leq b$, we have either
  $\rdn F(a,b)\geq(1-\delta)F(a,b)$ or $\rdn F(a,b) > \sqrt{n}$.  Thus, the cost
  of a tree with given level sizes is either within $(1-\delta)$ of the
  unmodified cost, or else is at least~$\sqrt{n}$.  Since the optimal
  choice is smaller than $\sqrt{n}$ (here we use $k>1$), the proof of
  the lower bound carries over unchanged.
  
  We consider the graph of light edges, which is $G_{n,p}$ with
  $p=\eps+o(\eps)$.  If every set~$A$ of size $|A|=a$ has at least
  $b$ neighbors in the light-edge graph, then $\rdn F(a,b)=F(a,b)$.
  (To see this, consider the sets $A$ and $B$ for which $|A|=a$,
  $|B|=b$, and $\rdn F(A,B)=\rdn F(a,b)$. If there were a heavy edge
  from $A$ to $B$, then we could delete the endpoint of that edge
  from~$B$, and replace it with a vertex not already in $B$ which is
  connected to $A$ via a light edge, and $\rdn F(A,\text{modified
    $B$})<\rdn F(A,B)$, a contradiction. Hence there is no heavy edge
  from $A$ to $B$, so $F(a,b)\leq F(A,B)=\rdn F(A,B)=\rdn F(a,b)\leq
  F(a,b)$.)

  We consider several cases.

  \paragraph{Case $b\leq \eps n/4$:}
  For any vertex, its degree $D$ in the light-edge graph is a binomial
  distribution with parameters $n-1$ and $p=(1+o(1))\eps$.  Since $\E[D]=(1+o(1)) n \eps$, the standard
  large-deviation formula (see \cite[Eqn.~5.6]{MR1036755}) that we used
  earlier tells us that
  $\Pr[D<\eps n/2]\leq e^{-(1+o(1))\E[D]/8} = e^{-n\eps/(8+o(1))}$.  A union
  bound then tells us that w.h.p.\ the minimal degree is at least
  $\eps n/2$.  Conditional on this event, any set~$A$ has at least
  $\eps n/2 - |A|$ neighbors, so if $a\leq b \leq \eps n/4$, it follows
  that $\rdn F(a,b)=F(a,b)$.

  \paragraph{Case $a\geq n^{1/3}$, $b\leq n-n^{3/4}$:}
  We argue that any disjoint sets $A$ and $C$ of sizes at least
  $n^{1/3}$ and $n^{3/4}$ have an edge between them.  This is a union
  bound over all pairs of sets: the number of pairs of sets is at most
  $3^n$, but each pair has no edge with probability
  $(1-p)^{n^{13/12}}$.  This implies that, for $n$ large enough, w.h.p.\ 
  $\rdn F(a,b)=F(a,b)$ for any
  $a\geq n^{1/3}$ and $b\leq n-n^{3/4}$, since any such set~$A$ has at
  most $n^{3/4}$ non-neighboring vertices in the light-edge graph.

  \paragraph{Case $a\geq n^{1/3}$, $b> n-n^{3/4}$:}
  By the monotonicity of $F$ and $\rdn F$, and the above case, for large
  enough~$n$ we have w.h.p.
  \[
  F(a,b) \geq \rdn F(a,b) \geq \rdn F(a,n-n^{3/4}) = F(a,n-n^{3/4}).
  \]
  However, based on our bounds on $F$ from
  Lemmas~\ref{L:no_cheap_sets} and \ref{L:no_cheap_sets2}, we have
  \[
  F(a,n-n^{3/4}) \geq (1-\delta)F(a,b)
  \]
  for any $\delta$ given $n$ large enough.

  \paragraph{Case $a\leq n^{1/3}$, $b>\eps n/4$:}
  By the monotonicity of $F$ and $\rdn F$, and the second case above,
  for large enough~$n$ we have w.h.p.
  \[
  \rdn F(a,b) \geq \rdn F(n^{1/3},\eps n/4) = F(n^{1/3},\eps n/4).
  \]
  By \lref{L:no_cheap_sets}, w.h.p.\ this is at least
  \[
  (1-\delta) \frac{\eps^2 n^2/16}{2 n n^{1/3}} \gg \sqrt{n} \gg \frac32
  n^{1/3},
  \]
  i.e., it exceeds the weight of the greedy spanning tree.
\end{proof}

\section{Open problems}

We identified a sharp cutoff of depth $\log_2\log n \pm \Theta(1)$
above which the minimum bounded-depth spanning tree has weight that is
asymptotically equal to the value of the unconstrained minimum
spanning tree, and below which it is much larger.  This same cutoff at
$\log_2\log n \pm \Theta(1)$ holds for minimum bounded-depth Steiner
trees with $m$ terminals when $m=\Theta(n)$, but we do not know the location
of the cutoff (or indeed if there is one) when $m$ is much smaller than $n$.
If there is a cutoff, we know that it occurs when the depth~$k$ 
is in the interval
$$\log_2\log m - \log_2\log(e n/m) - \omega(1) \leq k \leq 
 \log_2\log n + \omega\left(\frac{n}{m \log(e n/m)}\right),
$$
but we do not know where in the interval.  It would be interesting to
better understand the weights of bounded-depth Steiner trees for these
parameter values.

It would be interesting to understand better the large-$n$ behavior of
the weight of the bounded-depth MST near depth $\log_2 \log n +
\Delta$ as a function of $\Delta$.  The precise behavior could be
complicated, and is perhaps a periodic function of the fractional part
of $\log_2\log n$, but there are more basic open problems.  For
example, our construction in \sref{hybrid} shows that when the depth
bound is $\log_2\log n +\Delta$, the bounded-depth MST has weight
$\leq \zeta(3)+O(1/\sqrt{\Delta})$, while our best lower bound is
$\zeta(3)$.  We do not know how fast the approach to $\zeta(3)$ is
when $\Delta$ is increased, or indeed, if $\zeta(3)$ is reached for
some finite $\Delta$.

The weight of the minimum weight Steiner tree (with unbounded depth),
as a function of $\alpha=m/n$ (the ratio of the number of terminals to the
number of vertices) goes from~$0$ at $\alpha=0$ to $\zeta(3)$ at
$\alpha=1$.  As mentioned in \cite{BGRS}, it would be interesting to
understand how the weight varies from $0$ to $\zeta(3)$ for
intermediate values of~$\alpha$.

There was an experimental study aimed at sharpening our
estimate of $(1-o(1)) \frac32 n^{1/3}$ for the asymptotic weight of
minimum bounded-depth spanning trees when with depth bound $k=2$
\cite{BayatiBBCRZ2008}, suggesting $\frac32 n^{1/3} - \text{const}$.
This constant will depend on the weight distribution;
it may be interesting to rigorously determine the constant.

\section*{Acknowledgements}

This problem was called to our attention by Riccardo Zecchina and his
experimentation with others on the cavity method for bounded-depth
spanning trees and Steiner trees on random graphs.  Claire Mathieu
informed us of some relevant literature.

\pdfbookmark[1]{References}{bib}
\bibliographystyle{halpha}
\bibliography{bdmst}

\end{document}